%% file: FPT_1828_F.tex
\documentclass[10pt,twoside,reqno,draft]{amsart}
\input{prembu}
\usepackage{amsmath}
\usepackage{amssymb}
\usepackage{amsfonts}
\usepackage{amsthm}
\begin{document}
\setcounter{page}{1}

\leftline{\footnotesize {Accepted for publication in}}
\leftline{\footnotesize {\bf\em   Fixed Point Theory}}
\leftline{\footnotesize
http://www.math.ubbcluj.ro/$^{\sim}$nodeacj/sfptcj.html}

\vs*{1.5cm}

\title[Relation-theoretic Metrical Fixed Point Theorems]{\large Relation-theoretic metrical fixed point theorems under nonlinear contractions}
\author[Ahmadullah, Imdad and Gubran]{Md Ahmadullah$^{*}$, Mohammad Imdad, Rqeeb Gubran}
\date{04-09-2016}
\maketitle

\vs*{-0.5cm}

\bc
{\footnotesize
Department of Mathematics, Aligarh Muslim University,\\ Aligarh-202002, U.P., India.\\
E-mail: ahmadullah2201@gmail.com; mhimdad@gmail.com; rqeeeb@gmail.com\\
\medskip

}
\ec

\bigskip

{\footnotesize

\noindent
{\bf Abstract.}
We establish fixed point theorems for nonlinear contractions on a metric space (not essentially complete) endowed with an arbitrary binary relation. Our results extend, generalize, modify and unify several known results especially those contained in Samet and Turinici [Commun. Math. Anal. 13, 82-97 (2012)] and Alam and Imdad [J. Fixed Point Theory Appl. 17(4), 693-702 (2015)]. Interestingly a corollary to one of our main results proved under symmetric closure of any binary relation remains a sharpened version of a theorem due to Samet and Turinici. Finally, we use examples to highlight the realized improvements in the results proved in this paper.



}

\bigskip

\section{Introduction and Preliminaries}
Banach contraction principle \cite{Bnch1922} is one of the most fruitful and applicable theorems in classical functional analysis. In the last nine decades, this principle has been generalized and improved by numerous researchers in the different directions viz:
\begin{itemize}
  \item enlarging the class of underlying spaces,
  \item replacing contraction condition with relatively weaker contractive condition,
  \item weakening the involved metrical notions,
\end{itemize}
and such practice is still in business.

\bigskip
In 1986, the idea of order-theoretic fixed point results was initiated by
 Turinici \cite{Turinici1986, Turinicid1986}. In 2004, unknowingly,
Ran and Reunings \cite{RanR2004} rediscovered a slightly more natural order-theoretic version
 of Banach contraction principle. Recently, several authors utilized different types of
 binary relation viz: partial order (see Nieto and Rodr\'{\i}guez-L\'{o}pez \cite{NietoL2005}),
tolerance (see Turinici \cite{Turinici2011, Turinici2012}), strict order (see Ghods et al. \cite{GGGH2012}),
transitive (see Ben-El-Mechaiekh \cite{BenEM2015}), preorder (see Turinici \cite{Turinici2013}) $etc$
to prove their respective fixed point results.
Samet and Turinici \cite{SametT2012} established a fixed point theorem for
nonlinear contraction by using a symmetric closure of any binary relation.
 Most recently, Alam and Imdad \cite{Alamimdad,
Alamimdad2}  established a relation-theoretic version of Banach
contraction principle employing amorphous relation which in
turn generalizes several well known relevant order-theoretic fixed
point theorems. For more details, one can consult (\cite{AhmadJI,Alamimdad3,Alamimdad,Alamimdad2,BenEM2015,GGGH2012,Jach,
NietoL2005,RanR2004,SametT2012,Turinici1986, Turinicid1986,Turinici2013,Turinici2011,Turinici2012} and references cited therein).

\bigskip
\indent In the sequel, the following definitions will be utilized.

\noindent{\bf Definition 1.1} \cite{Lips1964}. A binary relation on a non-empty set $X$ is defined as a subset of $X\times X$ which will be denoted by $\mathcal{R}$. We say that ``$x$ is $\mathcal{R}$-related to $y$" if and only if $(x,y)\in \mathcal{R}$.

\bigskip In what follows, $\mathcal{R}$, $\mathbb{N}$ and $\mathbb{N}_{0}$ respectively, stand for a non-empty binary
relation, the set of natural numbers and the set of whole numbers.

\bigskip \noindent
{\bf Definition 1.2} \cite{Maddux2006}. A binary relation $\mathcal{R}$ defined on $X$ is called complete if for all $x,y\in X$, either $(x,y)\in
\mathcal{R}$ or $(y,x)\in \mathcal{R}$ which is denoted as $[x,y]\in \mathcal{R}$.

\bigskip \noindent
{\bf Definition 1.3} \cite{Lips1964}. Let $\mathcal{R}$ be a binary relation defined on a non-empty set $X$. Then the symmetric closure of $\mathcal{R}$ is defined as the smallest symmetric relation containing $\mathcal{R}$ ($i.e., \mathcal{S}:=\mathcal{R}\cup \mathcal{R}^{-1}$). Often, it is denoted by $\mathcal{S}$ or $\mathcal{R}^s$.

\bigskip \noindent
{\bf Definition 1.4} \cite{Alamimdad}. Let $f$ be a self-mapping defined on a non-empty
set $X$. Then a binary relation $\mathcal{R}$
defined on $X$ is called $f$-closed if
$$(x,y)\in \mathcal{R}\Rightarrow (fx,fy)\in \mathcal{R},~ {\rm for~ all}~ x,y\in X.$$

\bigskip \noindent
{\bf Definition 1.5} \cite{Alamimdad}. Let $\mathcal{R}$ be a binary relation defined on a non-empty set $X$.
 Then a sequence $\{x_n\}$ in $X$ is called $\mathcal{R}$-preserving if
$$(x_n,x_{n+1})\in\mathcal{R},\;\;{\rm for~all}~n\in \mathbb{N}.$$

\bigskip \noindent
{\bf Definition 1.6} \cite{Alamimdad2}. Let $(X,d)$ be a metric space and $\mathcal{R}$ a binary relation defined on $X$.
 We say that $(X,d)$ is $\mathcal{R}$-complete if every $\mathcal{R}$-preserving Cauchy sequence in $X$ converges to a point in $X$.

 \noindent
{\bf Remark 1.7} \label{rem2}\cite{Alamimdad2}. {\it Every complete metric space is
$\mathcal{R}$-complete, where $\mathcal{R}$ denotes a binary relation.
Particularly, if $\mathcal{R}$ is universal relation, then notions of
completeness and $\mathcal{R}$-completeness coincide.}

\bigskip \noindent
{\bf Definition 1.8} \cite{Alamimdad2}. Let $(X,d)$ be a metric space and $\mathcal{R}$
a binary relation defined on $X$. Then a mapping
$f:X\rightarrow X$ is called $\mathcal{R}$-continuous at $x$ if for
any $\mathcal{R}$-preserving sequence $\{x_n\}$ with
$x_n\stackrel{d}{\longrightarrow} x$, we have
$f(x_n)\stackrel{d}{\longrightarrow} f(x)$. Moreover, $f$ is called
$\mathcal{R}$-continuous if it is $\mathcal{R}$-continuous at every point of $X$.

 \noindent
{\bf Remark 1.9} \label{rem3} \cite{Alamimdad2}. {\it Every continuous mapping is
$\mathcal{R}$-continuous, where $\mathcal{R}$ denotes a binary relation.
Particularly, if $\mathcal{R}$ is universal relation, then notions of
continuity and $\mathcal{R}$-continuity coincide.}

\bigskip \noindent
{\bf Definition 1.10} \cite{SametT2012}. Let $(X,d)$ be a metric space and $\mathcal{S}$
 be the symmetric closure of a binary relation $\mathcal{R}$
   defined on $X$. We say that $(X,d,\mathcal{S})$ is regular if for any sequence
$\{x_n\}$ with $(x_n,x_{n+1})\in \mathcal{S}$ (for all $n\in\mathbb{N}$) and  $x_n\stackrel{d}{\longrightarrow} x$, then there
is a subsequence $\{x_{n_k}\}{\rm \;of\;} \{x_n\}$
such that $(x_{n_k},x)\in\mathcal{S},~~~{\rm for~all}~k\in \mathbb{N}.$

\bigskip \noindent
{\bf Definition 1.11} \cite{Alamimdad}. Let $(X,d)$ be a metric space and $\mathcal{R}$ a binary relation defined on $X$. Then $\mathcal{R}$ is called
$d$-self-closed if for any $\mathcal{R}$-preserving sequence
$\{x_n\}$ with $x_n\stackrel{d}{\longrightarrow} x$, there
is a subsequence $\{x_{n_k}\}{\rm \;of\;} \{x_n\} \;{\rm
such ~that}\;\;[x_{n_k},x]\in\mathcal{R},~~~{\rm for ~all}~k\in \mathbb{N}.$
\noindent Notice that, $(X,d,\mathcal{S})$ is regular if and only if $\mathcal{S}$ is
$d$-self-closed.

\bigskip \noindent
{\bf Definition 1.12} \cite{SametT2012}. Let $(X,d)$ be a metric space and $\mathcal{R}$ a binary relation
   defined on $X$. Then a subset $D$ of $X$ is called $\mathcal{R}$-directed if for every pair of points $x,y \in D$, there
is $z \in X$ such that $(x,z)\in\mathcal{R}$ and $(y,z)\in\mathcal{R}$.

\bigskip \noindent
{\bf Definition 1.13} \cite{KBR2000}. Let $(X,d)$ be a metric space, $\mathcal{R}$ a
binary relation defined on $X$ and $x,y$ a pair of points in $X$. Then a finite sequence
$\{z_0,z_1,z_2,...,z_{l}\}$ in $X$ is said to be a path
of length $l$ (where $l\in \mathbb{N}$) joining
$x$ to $y$ in $\mathcal{R}$ if $z_0=x, z_l=y$ and $(z_i,z_{i+1})\in\mathcal{R}$ for each $i\in \{1,2,3,\cdots ,l-1\}.$
\noindent Observe that, a path of length $l$ involves $(l+1)$ elements of $X$
that need not be distinct in general.

\bigskip \indent Given a metric space $(X,d)$, a self-mapping $f$ on $X$ and a binary relation $\mathcal{R}$ on $X$, we employ the following notations:
\begin{itemize}
  \item $F(f)$: the set of all fixed points of $f$;
  \item $X(f,\mathcal{R})$: the collection of all points $x \in X$ such that $(x,fx)\in \mathcal{R}$;
  \item $\Upsilon(x,y,\mathcal{R})$: the family of all paths joining $x$ to $y$ in $\mathcal{R}$;
  \item $\mathcal{M}_{f}(x,y)$ := $max\big\{d(x,y), d(x,fx), d(y,fy), \frac{1}{2}[d(x,fy)+d(y,fx)]\big\}$; and
  \item $\mathcal{N}_{f}(x,y)$ := $max\big\{d(x,y), \frac{1}{2}[d(x,fx)+d(y,fy)], \frac{1}{2}[d(x,fy)+d(y,fx)]\big\}$.
\end{itemize}

\noindent {\bf Remark 1.14}. \label{rem1}
{\it Observe that,$~~{\rm for~all}~x,y\in X,~~ \mathcal{N}_{f}(x,y)\leq\mathcal{M}_{f}(x,y).$}

\bigskip \noindent Let $\Phi$ be the family of all mappings $\varphi : [0,\infty) \to [0,\infty)$ satisfying the following properties
\begin{description}
  \item[$(\Phi_1)$] $\varphi$ is increasing;
  \item[$(\Phi_2)$] $\displaystyle\sum_{n=1}^{\infty}\varphi^{n}(t)<\infty$ for each $t > 0$, where $\varphi^{n}$ is the $n$-th iterate of $\varphi$.
\end{description}

 Recall that, the necessary condition of any real convergent series $ \sum_{n}^{ }a_n$ is that $\displaystyle\lim_{n\to \infty}a_n=0$

\bigskip \indent
The following two lemmas are required in our subsequent discussion.

\noindent {\bf Lemma 1.15}\label{lem1} \cite{SametT2012}. {\it Let $\varphi\in\Phi$. Then for all $t>0,$ we have $\varphi(t)<t.$}

\bigskip \noindent
{\bf Lemma 1.16}.\label{lem2} {\it Let $(x,d)$ be a metric space and $f$ a self-mapping on $X$. Then for each $x \in X$,
$$\mathcal{M}_{f}(x,fx)\leq max \big\{d(x,f x),d(fx,f^2x)\big\}.$$}
{\bf Proof.} Let $x$ be an arbitrary element of $X$. Then
  \begin{eqnarray}
  \nonumber \mathcal{M}_{f}(x,fx) &=& max\big\{d(x,fx), d(x,fx), d(fx,f^2x), \frac{1}{2}[d(x,f^2x)+d(fx,fx)]\big\} \\
    \nonumber &\leq&  max\big\{d(x,fx), d(fx,f^2x), \frac{1}{2}[d(x,fx)+d(fx,f^2x)]\big\}\\
   \nonumber &\leq&  max\big\{d(x,fx), d(fx,f^2x), max\{d(x,fx),d(fx,f^2x)\}\big\}\\
   \nonumber &=& max\big\{d(x,fx),d(fx,f^2x)\big\}.~~\square
  \end{eqnarray}

\bigskip
For the sake of completeness, we record the following known relevant results:

\bigskip
\noindent{\bf Theorem 1.17} {\rm(Theorem 2.1, Samet and Turinici \cite{SametT2012})}. \label{A}
 {\it Let $(X,d)$ be a metric space, $\mathcal{R}$ a binary relation defined on $X$ with a symmetric closure $\mathcal{S}:=\mathcal{R}\cup\mathcal{R}^{-1}$ and $f$ a self-mapping on
$X$. Assume that the following conditions hold:\\
\indent\hspace{0.5mm}$(i)$ $(X,d)$ is complete;\\
\indent\hspace{0.5mm}$(ii)$ there exists $x_0 \in X$ such that $(x_0,fx_0)\in \mathcal{S}$;\\
\indent\hspace{0.5mm}$(iii)$ $\mathcal{S}$ is $f$-closed;\\
\indent\hspace{0.5mm}$(iv)$ $(X,d,\mathcal{S})$ is regular;\\
\indent\hspace{0.5mm}$(v)$ there exists $\varphi\in\Phi $ such that
$$d(fx,fy)\leq \varphi( \mathcal{N}_{f}(x,y)),$$
\indent\hspace{0.5mm}for all $x,y\in X\;\textrm{with}\; (x,y)\in \mathcal{S}.$\\
Then $f$ has a fixed point. Moreover, if in addition, $F(f)$ is $\mathcal{S}$-directed, then $f$ has a unique fixed point.}

\bigskip \noindent{\bf Theorem 1.18} {\rm (Alam and Imdad \cite{Alamimdad, Alamimdad2})}.\label{B}
   {\it Let $(X,d)$ be a metric space,
$\mathcal{R}$ a binary relation on $X$ and $f$ a self-mapping on
$X$. Suppose that the following
conditions hold:\\
\indent\hspace{0.5mm}$(i)$ there exists $Y\subseteq X, fX\subseteq Y\subseteq X$ such that $(Y,d)$ is $\mathcal{R}$-complete;\\
\indent\hspace{0.5mm}$(ii)$ $X(f,\mathcal{R})$ is non-empty;\\
\indent\hspace{0.5mm}$(iii)$ $\mathcal{R}$ is $f$-closed;\\
\indent\hspace{0.5mm}$(iv)$ either $f$ is $\mathcal{R}$-continuous or $\mathcal{R}|_{Y}$ is $d$-self-closed;\\
\indent\hspace{0.5mm}$(v)$ there exists $\alpha\in [0,1)$ such that \\
\indent\hspace{2.5cm}$d(fx,fy)\leq\alpha d(x,y)\;\;{\rm for~all}~ x,y\in X$ with $(x,y)\in \mathcal{R}$.\\
Then $f$ has a fixed point. Moreover, if\\
\indent\hspace{0.5mm} $(vi)$ $\Upsilon(x,y,\mathcal{R}^s)$ is
non-empty, for each $x,y ~ in~ X$,\\
then $f$ has a unique fixed point.}

\bigskip \noindent{\bf Proposition 1.19}\label{p1} \cite{Alamimdad}. {\it If $\mathcal{R}$ is a binary relation
defined on a non-empty set $X$, then (for all $x,y$ in X)
$$(x,y)\in\mathcal{R}^s\Longleftrightarrow [x,y]\in\mathcal{R}.$$}

 \noindent{\bf Proposition 1.20}.\label{p2}
  {\it Let $(X,d)$ be a metric space equipped with a binary relation $\mathcal{R}$ defined on $X$, $f$
   a self-mapping on $X$ and $\varphi\in\Phi $. Then the following conditions are equivalent:
   \begin{description}
     \item[(A)] $d(fx,fy)\leq\varphi(\mathcal{M}_{f}(x,y)) ~with ~(x,y)\in \mathcal{R},$
     \item[(B)] $d(fx,fy)\leq\varphi(\mathcal{M}_{f}(x,y)) ~with ~[x,y]\in \mathcal{R}.$
   \end{description}}

\noindent{\bf Proof.}
  The implication $(B)\Rightarrow (A)$ is straightforward.

\noindent To show that $(A)\Rightarrow (B)$,
  choose $x,y\in X$ such that $[x,y]\in \mathcal{R}$. If $(x,y)\in \mathcal{R}$, then $(B)$ immediately follows from $(A)$. Otherwise, if $(y,x)\in \mathcal{R}$, then by $(A)$ and the symmetry of the metric $d$, we have
\begin{eqnarray*}
   d(fx,fy)=d(fy,fx)&\leq& \varphi\big(max\big\{d(y,x), d(y,fy), d(x,fx), \frac{1}{2}[d(y,fx)+d(x,fy)]\big\}\big) \\
   &=& \varphi(max\big\{d(x,y), d(x,fx), d(y,fy), \frac{1}{2}[d(x,fy)+d(y,fx)]\big\}) \\
   &=& \varphi(\mathcal{M}_{f}(x,y)).
\end{eqnarray*}
  Hence, $(A)\Rightarrow (B).$ $\square$

\bigskip
The main results of this paper are based on the following motivations and observations:
\begin{enumerate}
  \item the main result of Samet and Turinici \cite{SametT2012} is improved by replacing the symmetric closure $\mathcal{S}$ of any binary relation with amorphous binary relation $\mathcal{R}$,
  \item the completeness of the whole space $X$ involved in Theorem 1.17 (due to Samet and Turinici \cite{SametT2012}) is replaced by relatively weaker notion namely: $\mathcal{R}$-completeness of any subspace $Y\subseteq X$, such that $fX\subseteq Y\subseteq X$. Observe that the completeness of the whole space is not needed,
  \item the regularity of $X$ involved in Theorem 1.17, is replaced by $d$-self-closedness of $\mathcal{R}$,
  \item the contraction conditions involved in Theorems 1.17 and 1.18  are replaced by relatively weaker nonlinear generalized \'{C}iri\'{c} contraction,
  \item a corollary to our main result deduced for the symmetric closure $\mathcal{S}$ of any binary relation remains a sharper version of  Theorem 1.17,
  \item some examples are furnished to highlight the realized improvement in the results of this paper.
\end{enumerate}

\bigskip

\section{Main results}
Now, we are equipped to prove our main result as follows:

\noindent
{\bf Theorem 2.1}. {\it Let $(X,d)$ be a metric space equipped with a binary relation $\mathcal{R}$ defined on $X$ and $f$ a self-mapping on
$X$. Assume that the following conditions hold:\\
\indent\hspace{0.5mm}$(i)$ there exists $Y\subseteq X, fX\subseteq Y\subseteq X$ such that  $(Y,d)$ is $\mathcal{R}$-complete;\\
\indent\hspace{0.5mm}$(ii)$ $X(f,\mathcal{R})$ is non-empty;\\
\indent\hspace{0.5mm}$(iii)$ $\mathcal{R}$ is $f$-closed;\\
\indent\hspace{0.5mm}$(iv)$ either $f$ is $\mathcal{R}$-continuous or $\mathcal{R}|_Y$ is $d$-self-closed;\\
\indent\hspace{0.5mm}$(v)$ there exists $\varphi\in\Phi $ such that
$$d(fx,fy)\leq \varphi( \mathcal{M}_{f}(x,y)),$$
\indent\hspace{0.5mm}for all $x,y\in X\;\textrm{with}\; (x,y)\in \mathcal{R}.$\\
Then $f$ has a fixed point.}

\noindent
{\bf Proof.} Since $X(f,\mathcal{R})\neq \emptyset$. Let $x_0 \in X(f,\mathcal{R})$. Construct a Picard sequence $\{x_n\}$, with the initial point $x_{0},~ i.e.,$ $$x_{n+1}=f(x_{n}),\; \;{\rm for~all}~ n\in \mathbb{N}_0.\eqno(1)$$
Since $(x_0,fx_0)\in \mathcal{R}$ and $\mathcal{R}$ is $f$-closed, we have
$$(fx_0,f^2x_0),(f^2x_0,f^3x_0),\cdots ,(f^nx_0,f^{n+1}x_0),\cdots\in \mathcal{R}.$$
Thus,
$$(x_n,x_{n+1})\in \mathcal{R},\;\;{\rm for ~all}~n\in \mathbb{N}_0,$$
so that the sequence $\{x_n\}$ is $\mathcal{R}$-preserving.
From condition $(v)$, we have (for all $n\in \mathbb{N}$)
$$d(x_n,x_{n+1})=d(fx_{n-1},fx_{n})\leq \varphi( \mathcal{M}_{f}(x_{n-1},x_{n}))\eqno(2)$$
On using Lemma 1.16, we have (for all $n\in \mathbb{N}$)
$$\mathcal{M}_{f}(x_{n-1},x_{n})\leq max \big\{d(x_{n-1},x_{n}), d(x_n,x_{n+1})\big\}.\eqno(3)$$
On using (2), (3) and the property $(\Phi_1)$, we obtain (for all $n\in \mathbb{N}$)
$$d(x_n,x_{n+1})\leq \varphi\big( max \big\{d(x_{n-1},x_{n}), d(x_n,x_{n+1})\big\}\big).\eqno(4)$$
Now, we show that the sequence $\{x_n\}$ is Cauchy in $(X,d)$. In case $x_r=x_{r+1}$ for some $r\in \mathbb{N}_0,$ then the result is immediate. Otherwise, $x_n \neq x_{n+1}$ for all $n\in \mathbb{N}_0.$
Suppose that $d(x_{s-1},x_{s})\leq d(x_s,x_{s+1}),~ {\rm for ~some}~ s\in \mathbb{N}$. On using (4) and Lemma 1.15, we get
$$d(x_{s},x_{s+1})\leq\varphi( d(x_s,x_{s+1}) )< d(x_s,x_{s+1}),$$
which is a contradiction. Thus $d(x_n,x_{n+1})< d(x_{n-1},x_{n})$  (for all $n\in \mathbb{N}$), so that
$$d(x_{n},x_{n+1})\leq\varphi( d(x_{n-1},x_{n})),~ {\rm for ~all}~ n\in \mathbb{N}.$$
By induction on $n$ and the property $(\Phi_1)$, we get
$$d(x_{n},x_{n+1})\leq\varphi^{n}( d(x_0,x_1)),~ {\rm for ~all}~ n\in \mathbb{N}_0.$$
Now, for all $ m,n\in \mathbb{N}_0$ with $m\geq n,$ we have
\begin{eqnarray*}
 \nonumber d(x_{n},x_{m})&\leq& d(x_{n},x_{n+1})+d(x_{n+1},x_{n+2})+\cdots+d(x_{m-1},x_{m})\\
&\leq&\varphi^{n}(d(x_0,x_1))+\varphi^{n+1}(d(x_0,x_1))+\cdots+\varphi^{m-1}(d(x_0,x_1))\\
 &=& \sum\limits_{k=n}^{m-1} \varphi^{k}(d(x_0,x_1))\\
 &\leq& \sum\limits_{k\geq n} \varphi^{k}(d(x_0,x_1))\\
&\rightarrow& 0\;{\rm as}\; n\rightarrow \infty,
\end{eqnarray*}
which amounts to saying that the sequence $\{x_n\}$ is Cauchy in $X$. Hence, $\{x_n\}$ is $\mathcal{R}$-preserving Cauchy sequence in $X$.
Since, $\{x_n\}\subseteq fX\subseteq Y$ (due to (1) and $(i)$), therefore $\{x_n\}$ is $\mathcal{R}$-preserving Cauchy sequence in $Y$. As $(Y,d)$ is $\mathcal{R}$-complete, there exists $p\in Y$ such that $x_n\stackrel{d}{\longrightarrow} p$.

If $f$ is $\mathcal{R}$-continuous, then
$$p=\displaystyle\lim_{n\to\infty}x_{n+1}=\displaystyle\lim_{n\to\infty}fx_{n}=f\displaystyle\lim_{n\to\infty}x_{n}= fp.$$
Hence $p$ is the fixed point of $f$.

Alternatively, if $\mathcal{R}|_Y$ is $d$-self-closed, then for any $\mathcal{R}$-preserving sequence
$\{x_n\}$ in $Y$ with $x_n\stackrel{d}{\longrightarrow} p$, there is a subsequence
$\{x_{n_k}\}{\rm \;of\;} \{x_n\} \;{\rm
such~that}\;\;[x_{n_k},p]\in\mathcal{R}|_Y\subseteq \mathcal{R},~~~{\rm for~all}~ k\in \mathbb{N} _0.$

Write $\delta:=d(fp,p)\geq0$. Suppose on contrary that $\delta>0.$ On using condition $(v)$,
Proposition 1.20 and $[x_{n_k},p]\in\mathcal{R},$ for all $k\in \mathbb{N} _0,$ we get
$$d(x_{{n_k}+1},fp)=d(fx_{{n_k}},fp)\leq \varphi( \mathcal{M}_{f}(x_{{n_k}},p)),\eqno(5)$$
where $$\mathcal{M}_{f}(x_{{n_k}},p)= max\big\{d(x_{{n_k}},p), d(x_{{n_k}},x_{{n_k}+1}), d(p,fp), \frac{1}{2}[d(x_{{n_k}},fp) +d(p,x_{{n_k}+1})]\big\}.$$
If $\mathcal{M}_{f}(x_{{n_k}},p)=d(p,fp)=\delta,$ then (5) reduces to
$$d(x_{{n_k}+1},fp)\leq \varphi(\delta),$$
which on making $k \to \infty$, gives arise
$$\delta\leq \varphi(\delta),$$
which is a contradiction. Otherwise, if
$$\mathcal{M}_{f}(x_{{n_k}},p)= max\big\{d(x_{{n_k}},p), d(x_{{n_k}},x_{{n_k}+1}),
\frac{1}{2}[d(x_{{n_k}},fp) +d(p,x_{{n_k}+1})]\big\},$$ then due to the fact that $x_n\stackrel{d}{\longrightarrow} p$,
there exists $h=h(\delta)$ such that
$$\mathcal{M}_{f}(x_{{n_k}},p)\leq\frac{2}{3}\delta, ~{\rm {for~ all}}~k\geq h.$$
As $\varphi$ is increasing, we have
$$\varphi(\mathcal{M}_{f}(x_{{n_k}},p))\leq\varphi(\frac{2}{3}\delta), ~{\rm {for~ all}}~k\geq h.\eqno(6)$$
On using (5) and (6), we get
$$d(x_{{n_k}+1},fp)=d(fx_{{n_k}},fp)\leq\varphi(\frac{2}{3}\delta), ~{\rm {for~ all}}~k\geq h.$$
Letting $k \to \infty$ and using Lemma 1.15, we get
$$\delta\leq\varphi(\frac{2}{3}\delta)<\frac{2}{3}\delta<\delta,$$
which is again a contradiction. Hence, $\delta =0$, so that
 $$d(fp,p)=\delta =0\Rightarrow fp=p,$$
which concludes the proof.  ${\square}$

\bigskip

In particular, on setting $Y=X$ in Theorem 2.1, we deduce a corollary which is an improved version of Theorem 1.17 (up to fixed point) due to the involvement of relatively weaker notions in the considerations of completeness, regularity and contraction condition:

\noindent{\bf Corollary 2.2}.\label{cor0}
{\it Let $(X,d)$ be a metric space equipped with a binary relation $\mathcal{R}$ defined on $X$ and $f$ a self-mapping on
$X$. Assume that the following conditions hold:\\
\indent\hspace{0.5mm}$(i)$ $(X,d)$ is $\mathcal{R}$-complete;\\
\indent\hspace{0.5mm}$(ii)$ $X(f,\mathcal{R})$ is non-empty;\\
\indent\hspace{0.5mm}$(iii)$ $\mathcal{R}$ is $f$-closed;\\
\indent\hspace{0.5mm}$(iv)$ either $f$ is $\mathcal{R}$-continuous or $\mathcal{R}$ is $d$-self-closed;\\
\indent\hspace{0.5mm}$(v)$ there exists $\varphi\in\Phi $ such that
$$d(fx,fy)\leq \varphi( \mathcal{M}_{f}(x,y)),$$
\indent\hspace{0.5mm}for all $x,y\in X\;\textrm{with}\; (x,y)\in \mathcal{R}.$\\
Then $f$ has a fixed point.}

\bigskip
In view of Remarks 1.7 and 1.9, we deduce the following relatively more natural consequence of Theorem 2.1.

\noindent {\bf Corollary 2.3}.\label{cor1} {\it Let $(X,d)$ be a metric space equipped with
 a binary relation $\mathcal{R}$ defined on $X$ and $f$ a self-mapping on
$X$. Assume that the following conditions hold:\\
\indent\hspace{0.5mm}$(i)$ there exists $Y\subseteq X, fX\subseteq Y\subseteq X$ such that  $(Y,d)$ is complete;\\
\indent\hspace{0.5mm}$(ii)$ $X(f,\mathcal{R})$ is non-empty;\\
\indent\hspace{0.5mm}$(iii)$ $\mathcal{R}$ is $f$-closed;\\
\indent\hspace{0.5mm}$(iv)$ either $f$ is continuous or $\mathcal{R}|_Y$ is $d$-self-closed;\\
\indent\hspace{0.5mm}$(v)$ there exists $\varphi\in\Phi $ such that
$$d(fx,fy)\leq \varphi( \mathcal{M}_{f}(x,y)),$$
\indent\hspace{0.5mm}for all $x,y\in X\;\textrm{with}\; (x,y)\in \mathcal{R}.$\\
Then $f$ has a fixed point.}

\bigskip
In view of Remark 1.14, the following consequence of Theorem 2.1 is predictable.

\noindent {\bf Corollary 2.4}.\label{cor2}
{\it Let $(X,d)$ be a metric space equipped with a binary relation $\mathcal{R}$ defined on $X$ and $f$ a self-mapping on
$X$. Assume that the following conditions hold:\\
\indent\hspace{0.5mm}$(i)$ there exists $Y\subseteq X, fX\subseteq Y\subseteq X$ such that  $(Y,d)$ is $\mathcal{R}$-complete;\\
\indent\hspace{0.5mm}$(ii)$ $X(f,\mathcal{R})$ is non-empty;\\
\indent\hspace{0.5mm}$(iii)$ $\mathcal{R}$ is $f$-closed;\\
\indent\hspace{0.5mm}$(iv)$ either $f$ is $\mathcal{R}$-continuous or $\mathcal{R}|_Y$ is $d$-self-closed;\\
\indent\hspace{0.5mm}$(v)$ there exists $\varphi\in\Phi $ such that
$$d(fx,fy)\leq \varphi\big(max\big\{d(x,y), \frac{1}{2}[d(x,fx)+d(y,fy)], \frac{1}{2}[d(x,fy)+d(y,fx)]\big\}\big),$$
\indent\hspace{0.5mm}for all $x,y\in X\;\textrm{with}\; (x,y)\in \mathcal{R}.$\\
Then $f$ has a fixed point.}

\bigskip
Now, we prove the following results ensuring the uniqueness of the fixed point (corresponding to Corollary 2.4):

\noindent {\bf Theorem 2.5}.\label{th2}
 {\it In addition to the hypotheses of Corollary 2.4, suppose that the following condition holds:
 $$(vi):~ fX ~is~\mathcal{R}^s\text{-}directed.$$
Then $f$ has a unique fixed point.}

\noindent{\bf Proof.} By Corollary 2.4, $F(f)$ is non-empty. If $F(f)$ is singleton, then there is nothing to prove. Otherwise, to accomplish the proof, take two arbitrary elements $p,q ~{\rm in}~ F(f),$ so that
$$fp=p\;{\rm and}\;fq=q.$$
Now, we are required to show that $p=q$.
Since $F(f)\subseteq fX$ and $fX~{\rm is}~\mathcal{R}^s$-directed, therefore there exists $z\in X$ such that $[p,z]\in\mathcal{R}$ and $[q,z]\in\mathcal{R}$. Now, we construct a Picard  sequence $\{z_n\}$ corresponding to $z_0=z$, so that $z_{n}=f^{n}z_{0} ~{\rm for~all}~ n\in\mathbb{N}_0$.

Our claim is that $\displaystyle\lim_{n\to\infty}d(p,z_n)=0.$ If $d(p,z_s)=0,$ for some $s\in\mathbb{N}_0$, then the result is immediate. Otherwise, suppose $d(p,z_n)>0,$ for all $n\in\mathbb{N}_0.$
\noindent Since $[p, z_{n}]\in \mathcal{R}$, for all $n\in \mathbb{N}_0$ (due to the fact that $\mathcal{R}$ is $f$-closed and $[p, z]\in \mathcal{R}$), therefore on using Proposition 1.20 and hypothesis $(v)$, we have
$$d(p, z_{{n}+1})=d(fp, fz_{{n}})\leq \varphi( \mathcal{N}_{f}(p, z_{{n}})),\eqno(7)$$
where, \begin{eqnarray*}
\mathcal{N}_{f}(p, z_{{n}}) &=& max\big\{ d(p,z_{{n}}), \frac{1}{2}[d(p,p)+d(z_{{n}},z_{{n}+1})],
\frac{1}{2}[d(p,z_{{n}}) +d(p,z_{{n}+1})]\big\} \\
 &\leq& max\big\{d(p,z_n), \frac{1}{2}[d(p,z_{{n}}) +d(p,z_{{n}+1})]\big\}\\
 &\leq&max\big\{d(p,z_{{n}}), d(p,z_{{n}+1})\big\}.
\end{eqnarray*}
Using the property $(\Phi_1)$ and (7) (for all $n\in\mathbb{N}_0$), we get
\begin{eqnarray*}
  d(p, z_{{n}+1}) &\leq& \varphi\big( max\big\{d(p,z_{{n}}), d(p,z_{{n}+1})\big\}\big)\\
   &= & \varphi\big(d(p,z_{{n}})\big),
\end{eqnarray*}
otherwise, the fact that $\varphi\big(d(p,z_{{n+1}})\big)< d(p,z_{{n+1}})$ is contradicted.
So, by induction on $n$ and increasingness of $\varphi$, we get
$$d(p,z_{{n}})\leq\varphi^{n} d(p,z_{{0}}), {\rm ~for~ all} ~n\in \mathbb{N}_0,$$
which on making $n\to \infty $ and using the property $(\Phi_2)$, we get
$$\displaystyle\lim_{n\to\infty}d(p,z_n)=0.\eqno(8)$$
Similarly, we can prove that
$$\displaystyle\lim_{n\to\infty}d(q,z_n)=0.\eqno(9)$$
Using (8) and (9), we have
\begin{eqnarray*}
  d(p,q) &\leq& d(p,z_n)+d(z_n, q)  \\
   &\to& 0,~ \text{as} ~n\to \infty\\
   \Rightarrow p &=& q.
\end{eqnarray*}
Thus, $f$ has a unique fixed point. $\square$

\bigskip\noindent {\bf Remark 2.6}.\label{rem4}
{\it In Theorem 2.5, we have used a relatively more natural condition ``$fX$ is $\mathcal{R}^s$-directed" instead of ``$F(f)~is~ \mathcal{R}^s$-directed" which is too restrictive. Our proof continue to hold even if we take ``$F(f)~is~ \mathcal{R}^s$-directed". Sometimes it is difficult to find $F(f).$}

\bigskip\noindent {\bf Theorem 2.7}.\label{th3}
  {\it Theorem 2.5 remains true, if we replace the condition $(vi)$ by the following condition:
 $$(vi)^{\prime}:\mathcal{R}|_{fX}~is~complete.$$}

\noindent {\bf Proof.} From Corollary 2.4, $F(f)\ne\emptyset$. If $F(f)$ is singleton, then there is nothing to prove. Otherwise, take two arbitrary but distinct elements $p,q ~{\rm in}~ F(f),$ so that
$$fp=p\;{\rm and}\;fq=q.$$
Since $\mathcal{R}|_{fX}$ is complete, therefore $[p,q]\in \mathcal{R}$. Using Proposition 1.20 and condition $(v)$, we get
\begin{eqnarray*}
  d(p,q)=d(fp,fq)&\leq& \varphi\big(max\big\{d(p,q), \frac{1}{2}[d(p,fp)+d(q,fq)], \frac{1}{2}[d(p,fq)+d(q,fp)]\big\}\big)\\
   &=& \varphi\big(max\big\{d(p,q), \frac{1}{2}[d(p,p)+d(q,q)], \frac{1}{2}[d(p,q)+d(q,p)]\big\}\big)\\
   &=& \varphi\big(d(p,q)\big).
\end{eqnarray*}
which contradicts to the fact that $\varphi\big(d(p,q)\big)< d(p,q)$ (in view of Lemma 1.15).
Hence $p=q$.
Thus $f$ has a unique fixed point. $\square$

\bigskip
In view of Remark 2.6, on choosing $\mathcal{R}$ to be the symmetric closure $\mathcal{S}$ of any arbitrary binary relation in Theorem 2.5, we obtain the following sharpened version of Theorem 1.17.

\noindent {\bf Corollary 2.8}.\label{cor3}
{\it Let $(X,d)$ be a metric space equipped with the symmetric closure $\mathcal{S}$ of any arbitrary binary relation defined on $X$ and $f$ a self-mapping on
$X$. Assume that the following conditions hold:\\
\indent\hspace{0.5mm}$(i)$ there exists $Y\subseteq X, fX\subseteq Y\subseteq X$ such that  $(Y,d)$ is $\mathcal{S}$-complete;\\
\indent\hspace{0.5mm}$(ii)$ $X(f,\mathcal{S})$ is non-empty;\\
\indent\hspace{0.5mm}$(iii)$ $\mathcal{S}$ is $f$-closed;\\
\indent\hspace{0.5mm}$(iv)$ either $f$ is $\mathcal{S}$-continuous or $\mathcal{S}|_Y$ is $d$-self-closed;\\
\indent\hspace{0.5mm}$(v)$ there exists $\varphi\in\Phi $ such that
$$d(fx,fy)\leq \varphi\big(max\big\{d(x,y), \frac{1}{2}[d(x,fx)+d(y,fy)], \frac{1}{2}[d(x,fy)+d(y,fx)]\big\}\big),$$
\indent\hspace{0.5mm}for all $x,y\in X\;\textrm{with}\; (x,y)\in \mathcal{S}.$\\
Then $f$ has a fixed point. Moreover, if in addition, $F(f)$ is $\mathcal{S}$-directed, then $f$ has unique fixed point.}

Notice that the hypothesis `$\mathcal{S}$ is $f$-closed' is equivalent to `the comparative property of $f$' and `$\mathcal{S}|_Y$ is $d$-self-closed' is equivalent to `the regular property of $(Y,d,\mathcal{S})$'.

\bigskip
Let $\Lambda$ be the collection of all increasing continuous mappings $\psi: [0,\infty)\to [0,\infty)$ such that
\begin{description}
  \item[{$(\Lambda_1)$}] For all $t>0,~~ 0<\psi(t)<t;$
  \item[$(\Lambda_2)$] $g(t)=\frac{t}{t-\psi(t)}$ is strictly decreasing function on $(0,\infty)$;
  \item[$(\Lambda_3)$] $\int_{0}^{T}g(t)dt< \infty $, for all $T>0$.
\end{description}

\noindent {\bf Lemma 2.9} \cite{Altman}.\label{lem3}
 {\it We have $\Lambda\subset \Phi$.}

\bigskip
In view of Lemma 2.9, we obtain the following consequence of Theorem 2.5:

\noindent {\bf Corollary 2.10}.
{\it The conclusion of Theorem 2.5 remains true if we replace the condition $(v)$ by the following (besides retaining rest of the hypotheses):
$$(v)^{\prime} ~there~ exists~ \psi\in\Lambda ~such ~that~ d(fx,fy)\leq \psi\big(\mathcal{N}_f(x,y)\big),$$
\indent\hspace{0.5mm}for all $x,y\in X\;\textrm{with}\; (x,y)\in \mathcal{R}$.}

Notice that, Corollary 2.10 is a sharpened version of Corollary 2.7 of Samet and Turinici \cite{SametT2012}.

\section{Consequences}
As consequences of our earlier established results, we derive several well known results of the existing literature.

\bigskip
\subsection{Relation-theoretic fixed point results}
On the setting $\varphi(t)=kt$, with $k\in [0,1)$, we derive the following corollaries which are immediate consequences of  Theorem 2.5.

\noindent{\bf Corollary 3.1}. {\it Let $(X,d)$ be a metric space equipped with a binary relation $\mathcal{R}$ defined on $X$ and $f$ a self-mapping on
$X$. Assume that the following conditions hold:\\
\indent\hspace{0.5mm}$(i)$ there exists $Y\subseteq X, fX\subseteq Y\subseteq X$ such that  $(Y,d)$ is $\mathcal{R}$-complete;\\
\indent\hspace{0.5mm}$(ii)$ $X(f,\mathcal{R})$ is non-empty;\\
\indent\hspace{0.5mm}$(iii)$ $\mathcal{R}$ is $f$-closed;\\
\indent\hspace{0.5mm}$(iv)$ either $f$ is $\mathcal{R}$-continuous or $\mathcal{R}|_Y$ is $d$-self-closed;\\
\indent\hspace{0.5mm}$(v)$ there exists $k\in [0,1)$ such that
$$d(fx,fy)\leq k\big(max\big\{d(x,y), \frac{1}{2}[d(x,fx)+d(y,fy)], \frac{1}{2}[d(x,fy)+d(y,fx)]\big\}\big),$$
\indent\hspace{0.5mm}for all $x,y\in X\;\textrm{with}\; (x,y)\in \mathcal{R}.$\\
Then $f$ has a fixed point. Moreover, if in addition, $fX$ is $\mathcal{R}^s$-directed, then $f$ has unique fixed point.}

\noindent {\bf Remark 3.2}. {\it Corollary 3.1 is a sharpened version of Corollary 2 (corresponding to condition (11)) due to Ahmadullah et al. \cite{AhmadJI}.}

\bigskip\noindent {\bf Corollary 3.3}. {\it Let $(X,d)$ be a metric space equipped with a binary relation $\mathcal{R}$ defined on $X$ and $f$ a self-mapping on
$X$. Assume that the following conditions hold:\\
\indent\hspace{0.5mm}$(i)$ there exists $Y\subseteq X, fX\subseteq Y\subseteq X$ such that  $(Y,d)$ is $\mathcal{R}$-complete;\\
\indent\hspace{0.5mm}$(ii)$ $X(f,\mathcal{R})$ is non-empty;\\
\indent\hspace{0.5mm}$(iii)$ $\mathcal{R}$ is $f$-closed;\\
\indent\hspace{0.5mm}$(iv)$ either $f$ is $\mathcal{R}$-continuous or $\mathcal{R}|_Y$ is $d$-self-closed;\\
\indent\hspace{0.5mm}$(v)$ there exist $a,b, c\geq 0$ with $a+2b+2c<1$ such that
$$d(fx,fy)\leq ad(x,y)+ b[d(x,fx)+d(y,fy)]+c[d(x,fy)+d(y,fx)],$$
\indent\hspace{0.5mm}for all $x,y\in X\;\textrm{with}\; (x,y)\in \mathcal{R}.$\\
Then $f$ has a fixed point. Moreover, if in addition, $fX$ is $\mathcal{R}^s$-directed, then $f$ has unique fixed point.}

\noindent {\bf Remark 3.4}. {\it Corollary 3.3 remains a sharpened version of a relation-theoretic \'{C}iri\'{c} fixed point theorem proved in Ahmadullah et al. \cite{AhmadJI} $viz.$ Corollary 2, corresponding to (13).}

\bigskip
The following result was obtained by Alam and Imdad \cite{Alamimdad2}:

\bigskip\noindent {\bf Corollary 3.5}.\label{cor7} {\it Let $(X,d)$ be a metric space equipped
 with a binary relation $\mathcal{R}$ defined on $X$ and $f$ a self-mapping on
$X$. Assume that the following conditions hold:\\
\indent\hspace{0.5mm}$(i)$ there exists $Y\subseteq X, fX\subseteq Y\subseteq X$ such that  $(Y,d)$ is $\mathcal{R}$-complete;\\
\indent\hspace{0.5mm}$(ii)$ $X(f,\mathcal{R})$ is non-empty;\\
\indent\hspace{0.5mm}$(iii)$ $\mathcal{R}$ is $f$-closed;\\
\indent\hspace{0.5mm}$(iv)$ either $f$ is $\mathcal{R}$-continuous or $\mathcal{R}|_Y$ is $d$-self-closed;\\
\indent\hspace{0.5mm}$(v)$ there exists $k\in [0,1)$ such that
$$d(fx,fy)\leq kd(x,y),$$
\indent\hspace{0.5mm}for all $x,y\in X\;\textrm{with}\; (x,y)\in \mathcal{R}.$\\
Then $f$ has a fixed point. Moreover, if in addition, $fX$ is $\mathcal{R}^s$-directed, then $f$ has unique fixed point.}

\bigskip\noindent {\bf Corollary 3.6}.\label{cor8} {\it Let $(X,d)$ be a metric space equipped with a binary relation $\mathcal{R}$ defined on $X$ and $f$ a self-mapping on
$X$. Assume that the following conditions hold:\\
\indent\hspace{0.5mm}$(i)$ there exists $Y\subseteq X, fX\subseteq Y\subseteq X$ such that  $(Y,d)$ is $\mathcal{R}$-complete;\\
\indent\hspace{0.5mm}$(ii)$ $X(f,\mathcal{R})$ is non-empty;\\
\indent\hspace{0.5mm}$(iii)$ $\mathcal{R}$ is $f$-closed;\\
\indent\hspace{0.5mm}$(iv)$ either $f$ is $\mathcal{R}$-continuous or $\mathcal{R}|_Y$ is $d$-self-closed;\\
\indent\hspace{0.5mm}$(v)$ there exists $k\in [0,1/2)$ such that
$$d(fx,fy)\leq k[d(x,fx)+d(y,fy)],$$
\indent\hspace{0.5mm}for all $x,y\in X\;\textrm{with}\; (x,y)\in \mathcal{R}.$\\
Then $f$ has a fixed point. Moreover, if in addition, $fX$ is $\mathcal{R}^s$-directed, then $f$ has unique fixed point.}

\noindent {\bf Remark 3.7}. {\it Corollary 3.6 remains an improved version of a relation-theoretic Kannan fixed point theorem ($i.e.,$ Corollary 2 corresponding to (9)) established in Ahmadullah et al. \cite{AhmadJI}.}

\bigskip\noindent {\bf Corollary 3.8}. \label{cor9} {\it Let $(X,d)$ be a metric space equipped with a binary relation $\mathcal{R}$ defined on $X$ and $f$ a self-mapping on
$X$. Assume that the following conditions hold:\\
\indent\hspace{0.5mm}$(i)$ there exists $Y\subseteq X, fX\subseteq Y\subseteq X$ such that  $(Y,d)$ is $\mathcal{R}$-complete;\\
\indent\hspace{0.5mm}$(ii)$ $X(f,\mathcal{R})$ is non-empty;\\
\indent\hspace{0.5mm}$(iii)$ $\mathcal{R}$ is $f$-closed;\\
\indent\hspace{0.5mm}$(iv)$ either $f$ is $\mathcal{R}$-continuous or $\mathcal{R}|_Y$ is $d$-self-closed;\\
\indent\hspace{0.5mm}$(v)$ there exists $k\in [0,1/2)$ such that
$$d(fx,fy)\leq k[d(x,fy)+d(y,fx)],$$
\indent\hspace{0.5mm}for all $x,y\in X\;\textrm{with}\; (x,y)\in \mathcal{R}.$\\
Then $f$ has a fixed point. Moreover, if in addition, $fX$ is $\mathcal{R}^s$-directed, then $f$ has unique fixed point.}

\noindent {\bf Remark 3.9} {\it Corollary 3.8 remains a sharpened version of a relation-theoretic Chatterjea fixed point theorem proved in Ahmadullah et al. \cite{AhmadJI}, Corollary 2 corresponding to (10).}

\bigskip
\subsection{Fixed point results in abstract space}

Under the universal relation ($i.e., \mathcal{R}=X\times X$), Theorems 2.5 and 2.7 unify to the following lone corollary:

\bigskip\noindent {\bf Corollary 3.10}.\label{cor10}  {\it Let $(X,d)$ be a metric space and $f$ a self-mapping on
$X$. Assume that the following conditions hold:\\
\indent\hspace{0.5mm}$(i)$ there exists $Y\subseteq X,$ $fX\subseteq Y\subseteq X$ such that $(Y,d)$ is complete;\\
\indent\hspace{0.5mm}$(iii)$ there exists $\varphi\in\Phi $ such that
$$d(fx,fy)\leq \varphi\big(max\big\{d(x,y), \frac{1}{2}[d(x,fx)+d(y,fy)], \frac{1}{2}[d(x,fy)+d(y,fx)]\big\}\big),$$
\indent\hspace{0.5mm}for all $x,y\in X.$\\
Then $f$ has a unique fixed point.}

Now, we consider the following special cases:
\begin{itemize}
  \item On setting $\varphi(t)=kt$ (with $k\in [0,1)$) in Corollary 3.10, we deduce an improved version of \'{C}iri\'{c} fixed point theorem \cite{Ciric1974}.
  \item Corollary 3.5 deduces to a sharpened form of Banach contraction principle \cite{Bnch1922}, under the universal relation.
  \item If $\mathcal{R}=X\times X$, Corollary 3.6 deduces to an improved version of Kannan fixed point theorem \cite{Kann1969}.
 \item Under the universal relation, Corollary 3.8 deduces to an improved version of Chatterjea fixed point theorem \cite{Chattj1972}.
\end{itemize}

\bigskip
\section{Illustrative Examples}
In this section, we furnish some examples to highlight the realized improvement in our results proved in this paper.

\bigskip\noindent {\bf Example 4.1.}\label{exm1} Consider a usual metric space $(X,d)$, where $X=(-1,4)$.
Now, define a binary relation $\mathcal{R}=\{(x,y)\in X^2~|~ x\geq y\}$, an increasing
mapping $\varphi : [0,\infty)\to [0,\infty)$ by $\varphi(t)= \frac{t}{2}$ and a self-mapping $f:X\to X$ by
$$f(x)=\left\{
  \begin{array}{ll}
    \frac{x}{2}, & \hbox{$x\in (-1,2]$;} \\
    1, & \hbox{$x\in (2,4)$.}
  \end{array}
\right.$$
Let $Y=[-\frac{1}{2},2)$, so that $fX=(-\frac{1}{2},1]\subset Y$ and $Y$ is $\mathcal{R}$-complete but $X$ is not $\mathcal{R}$-complete.
Evidently, $\mathcal{R}$ is $f$-closed and $f$ is $\mathcal{R}$-continuous.

By routine calculations, one can verify hypothesis $(v)$ of Theorems 2.5 and 2.7.
Moreover, as $fX ~is~\mathcal{R}^s\text{-}directed$ and $\mathcal{R}|_{fX}$ is complete, therefore all
 the hypotheses of Theorems 2.5 and 2.7 are satisfied, ensuring the uniqueness of the fixed point. Observe that, $x=0$ is the only point fixed of $f$.

With a view to establish genuineness of our results, notice that $\mathcal{R}=\{(x,y)\in X^2~|~ x\geq y\}$ is not symmetric and $\mathcal{R}$ can not  be a symmetric closure of any binary relation. Also $(X,d)$ is not complete and even not $\mathcal{R}$-complete,
which shows that the condition $(i)$ of Theorem 1.17 (due to Samet and Turinici \cite{SametT2012}) and even Corollary 2.2 is not satisfied.
Thus, our Theorems 2.1, 2.5 and 2.7 are applicable  to the present example
while Theorem 1.17 and even Corollary 2.2 are not, which substantiates the utility of Theorems 2.1, 2.5 and 2.7.

\bigskip\noindent {\bf Example 4.2.}\label{exm2}
  Consider $X=[0,3)$ equipped with usual metric, $i.e., d(x,y) = |x - y|$ for all $x,y\in X$ and a binary relation $\mathcal{R}=\{(0,0),(0,1),(0,2),(1,1),(2,2)\}$ on $X$, whose symmetric closure is $\mathcal{S}=\{(0,0),(0,1),(1,0),(0,2),(2,0),(1,1), (2,2)\}$. Define an increasing function $\varphi : [0,\infty)\to [0,\infty)$ by $\varphi(t)= \frac{3}{4}t$  and a self-mapping $f:X\to X$ by
$$f(x)=\left\{
  \begin{array}{ll}
     0, & \hbox{$x\in [0,1]$;} \\
    1, & \hbox{$x\in (1,3)$.}
  \end{array}
\right.$$

Let $Y=[0,1]$, so that $fX=\{0,1\}\subset Y$ and $Y$ is $\mathcal{S}$-complete.
Evidently, $f$ is not continuous but $\mathcal{S}$ is $f$-closed and $\varphi \in \Phi$. Take any $\mathcal{S}$-preserving sequence
$\{x_n\}$ in $Y, i.e.,$  $$(x_n,x_{n+1})\in\mathcal{S}, ~{\rm for~ all}~ n\in \mathbb{N}
~ {\rm with} ~x_n\stackrel{d}{\longrightarrow} x.$$
 Here, one may notice that $(x_n,x_{n+1})\in\mathcal{S}|_{Y}$, for all $n\in\mathbb{N}$
and there exists $N \in\mathbb{N}$
such that $x_n=x \in \{0,1\}, ~{\rm for~ all}~ n\geq N$. So, we need to choose a subsequence
$\{x_{n_k}\}$ of the sequence $\{x_n\}$ such that $x_{n_k}=x$, for all $k\in \mathbb{N}$, which
amounts to saying that $(x_{n_k},x)\in \mathcal{S}|_{Y}$, for all $k\in \mathbb{N}$. Therefore,
$\mathcal{S}|_{Y}$ is $d$-self-closed or $(Y,d,\mathcal{S})$ is regular.

Now, we check the condition $(v)$ of Theorem 1.17 and Corollary 2.8. For this, we need to verify for
$(x,y)\in \{(0,2),(2,0)\}.$ Otherwise, we have $d(fx,fy)=0,$ hence condition $(v)$ is obvious.
If $(x,y)\in \{(0,2),(2,0)\}$, then $d(fx,fy)=1 \leq \varphi(2)=\varphi(d(x,y))$ and
hence $f$ has a fixed point. As ${fX}$ is $\mathcal{S}$-directed, therefore all
the hypotheses of Corollary 2.8 are satisfied. Observe that, $x=0$ is the only fixed point of $f$.

Interestingly, Theorem 1.17 is not applicable to the present example as the underlying metric space $(X,d)$ is not complete, which shows that our result ($i.e.,$ Corollary 2.8) is an improvement over Theorem 1.17 (due to Samet and Turinici \cite{SametT2012}).

\bigskip\noindent {\bf Example 4.3.}\label{exm3}
 Let $X=[0,4]$ equipped with usual metric, $i.e.,$ for all $x,y\in X, d(x,y) = |x - y|$ and a binary relation $\mathcal{R}=\{(0,0),(1,1),(3,3),(4,4), (1,2),(3,4)\}$ on X. Define an increasing mapping $\varphi : [0,\infty)\to [0,\infty)$ by $\varphi(t)= \frac{3}{4}t$  and a self-mapping $f:X\to X$ by
$$f(x)=\left\{
  \begin{array}{ll}
    0, & \hbox{$x\in [0,1)$;} \\
    3, & \hbox{$x\in [1,2)$;}\\
    4, & \hbox{$x\in [2,4]$.}
  \end{array}
\right.$$

Clearly, $fX=\{0,3,4\}\subset Y=X$ where $Y$ is $\mathcal{R}$-complete, $f$ is not continuous but $\mathcal{R}$ is $f$-closed and $\varphi \in \Phi$.
On the lines of Example 4.2, one can verify that $\mathcal{R}|_{Y}$ is $d$-self-closed.

\vspace{.3cm}
Now, with a view to check the condition $(v)$ of Theorem 2.1, let $(x,y)=(1,2)$ (as in rest of cases $d(fx,fy)=0$), we have
$$d(f1,f2)=d(3,4)=1<\frac{3}{2} = \varphi(2)=\varphi\big(\frac{d(1,f2)+d(2,f1)}{2}\big),$$
which shows that condition $(v)$ (of Theorem 2.1) is satisfied,
and henceforth $f$ has a fixed point. Since $\{0,3\}\subset{fX}$ but $(0,3)\notin \mathcal{R}^s$, therefore $\mathcal{R}$ is not complete and hence ${fX}$ is not $\mathcal{R}^s$-directed. Observe that, $x=0, 4$ are the fixed points of $f$.

\vspace{.3cm} In the present example, observe that
$$(1,2)\in \mathcal{R} ~{\rm but}~ d(f1,f2) =d(3,4)\leq k d(1,2),~i.e.,~ 1\leq k$$
which shows that the contraction condition $(v)$ of Theorem 1.18, due to Alam and Imdad \cite{Alamimdad,Alamimdad2} is not satisfied. Also, $\mathcal{R}$ can not be a symmetric closure of any binary relation.

Thus, in all our results generalize, modify and unify the results of Samet and Turinici \cite{SametT2012} and Alam and Imdad
 \cite{Alamimdad, Alamimdad2}.

\bigskip
{\bf Acknowledgements:} All the authors are grateful to anonymous referee for his valuable suggestions and comments. Also, the first author is thankful to University Grant Commission, New Delhi, Govt. of India, for awarding Moulana Azad National Fellowship.

\bigskip

{\it Received: May 11, 2016; Accepted: September 01, 2016}

\end{document}

%% file: prembu.tex
\newcommand{\vs}{\vspace}
\newcommand{\bc}{\begin{center}}
\newcommand{\ec}{\end{center}}